\documentclass{article}
\usepackage{amsmath,amsfonts,amssymb,bm,hyperref}
\numberwithin{equation}{section}
\newtheorem{theorem}{Theorem}[section]

\newtheorem{remark}{Remark}[section]
\newtheorem{definition}{Definition}[section]

\usepackage{authblk}

\begin{document}

\title{Spectral theoretic characterization of the massless Dirac action}

\author{Robert J.~Downes\thanks{
Department of War Studies,
King's College London,
Strand,
London WC2R 2LS, UK;
robert.downes@kcl.ac.uk
}
}
\author{Dmitri Vassiliev
\thanks{
Department of Mathematics,
University College London,
Gower Street,
London WC1E~6BT, UK;
D.Vassiliev@ucl.ac.uk;
\url{http://www.homepages.ucl.ac.uk/\~ucahdva/};
supported by EPSRC grant EP/M000079/1
}
}
\affil{}

\maketitle

\begin{abstract}
We consider an elliptic self-adjoint first order differential
operator $L$ acting on pairs (2-columns) of complex-valued
half-densities over a connected compact 3-dimensional manifold
without boundary. The principal symbol of the operator $L$ is assumed to be
trace-free and the subprincipal symbol is assumed to be zero.
Given a positive scalar weight function, we study the
weighted eigenvalue problem for the operator $L$. The
corresponding counting function (number of eigenvalues between
zero and a positive $\lambda$) is known to admit,
under appropriate assumptions on periodic trajectories,
a two-term asymptotic expansion as $\lambda\to+\infty$ and we have recently
derived an explicit formula for the second asymptotic
coefficient. The purpose of this paper is to establish the
geometric meaning of the second asymptotic coefficient. To this
end, we identify the geometric objects encoded within our
eigenvalue problem --- metric, nonvanishing spinor field and
topological charge --- and express our asymptotic coefficients
in terms of these geometric objects. We prove that the second
asymptotic coefficient of the counting function has the
geometric meaning of the massless Dirac action.
\end{abstract}

\


\


\section{Main result}
\label{Main result}

Consider a first order differential
operator $L$ acting on 2-columns
$v=\begin{pmatrix}v_1&v_2\end{pmatrix}^T$
of complex-valued half-densities
over a connected compact 3-dimensional mani\-fold $M$ without boundary.
We assume the coefficients of the operator $L$ to be infinitely smooth. We also
assume that the operator $L$ is formally self-adjoint (symmetric):
$\int_Mu^*Lv\,dx=\int_M(Lu)^*v\,dx$ for all infinitely smooth
$u,v:M\to\mathbb{C}^2$. Here and further on
the superscript $\,{}^*\,$ in matrices, rows and columns
indicates Hermitian conjugation in $\mathbb{C}^2$
and $dx:=dx^1dx^2dx^3$, where $x=(x^1,x^2,x^3)$ are local
coordinates on $M$.

Let $L_\mathrm{prin}(x,p)$ be the principal symbol of the operator $L$,
i.e.~matrix obtained by leaving in $L$ only the leading (first order)
derivatives and replacing each $\partial/\partial x^\alpha$ by
$ip_\alpha$, $\alpha=1,2,3$. Here $p=(p_1,p_2,p_3)$ is the variable dual to the position
variable~$x$; in physics literature the $p$ would be referred to
as \emph{momentum}. Our principal symbol $L_\mathrm{prin}(x,p)$ is a $2\times 2$
Hermitian matrix-function on the cotangent bundle $T^*M$, linear in
every fibre $T_x^*M$ (i.e.~linear in $p$).
We assume that
$
\det L_\mathrm{prin}(x,p)\ne0
$
for all $(x,p)\in T'M:=T^*M\setminus\{p=0\}$
(cotangent bundle with the zero section removed),
which is a version of the ellipticity condition.

\begin{remark}
The tradition in microlocal analysis is to denote momentum by $\xi$.
We choose to denote it by $p$ instead because we will
need the letter $\xi$ for the spinor.
\end{remark}

We now make two additional assumptions:
\begin{itemize}
\item
we assume the principal symbol to be trace-free and
\item
we assume the subprincipal symbol of the operator $L$ to be zero
(see Appendix
\ref{Invariant analytic description of a first order differential operator}
for the definition of the subprincipal symbol).
\end{itemize}
The latter condition implies that our differential operator $L$
is completely determined by its principal symbol.
Namely, in local coordinates our operator reads
\[
L=-i[(L_\mathrm{prin})_{p_\alpha}(x)]\frac\partial{\partial x^\alpha}
-\frac i2(L_\mathrm{prin})_{x^\alpha p_\alpha}(x)\,,
\]
where the subscripts
indicate partial derivatives and
the repeated index $\alpha$ indicates summation over $\alpha=1,2,3$.
Of course, the above formula is a special case of formula
(\ref{operator in terms of its principal and subprincipal symbols}).

We study the eigenvalue problem
\begin{equation}
\label{eigenvalue problem}
Lv=\lambda wv,
\end{equation}
where $w(x)$ is a given infinitely smooth positive scalar weight function.
Obviously, the
problem (\ref{eigenvalue problem}) has the same spectrum as the problem
\begin{equation}
\label{eigenvalue problem without weight}
w^{-1/2}Lw^{-1/2}v=\lambda v,
\end{equation}
so it may appear that the weight function $w(x)$ is redundant.
We will, however, work with the eigenvalue problem (\ref{eigenvalue problem})
rather than with (\ref{eigenvalue problem without weight}) because we want our problem
to possess a gauge degree of freedom
(\ref{gauge transformation of eigenvalue problem}).
This gauge degree of freedom will eventually manifest itself as the conformal invariance of
the massless Dirac action, see Section~\ref{Conformal invariance} for details.

The problem (\ref{eigenvalue problem}) has a discrete spectrum accumulating to $\pm\infty$.
We define the counting function
$
N(\lambda):=\,\sum\limits_{0<\lambda_k<\lambda}1
$
as the number of eigenvalues $\lambda_k$
of the problem (\ref{eigenvalue problem}),
with account of multiplicities,
between zero and a positive~$\lambda$.
Theorem~8.4 from \cite{jst_part_a} states
that under appropriate assumptions on periodic trajectories
our counting function admits a two-term asymptotic expansion
\begin{equation}
\label{two-term asymptotic formula for counting function}
N(\lambda)=a\lambda^3+b\lambda^2+o(\lambda^2)
\end{equation}
as $\lambda\to+\infty$. If one wishes to reformulate
the asymptotic formula (\ref{two-term asymptotic formula for counting function})
in such a way that it remains valid without assumptions on periodic trajectories,
this can easily be achieved, say, by taking a convolution with a function from Schwartz space
$\mathcal{S}(\mathbb{R})$; see Theorem~7.2 in \cite{jst_part_a} for details.

Alternatively, one can look at the eta function
$\,\eta(s):=\sum|\lambda|^{-s}\operatorname{sign}\lambda\,$,
where summation is carried out over all nonzero eigenvalues
$\lambda$ and $s\in\mathbb{C}$ is the independent variable.
The series converges absolutely for
$\operatorname{Re}s>3$ and defines a holomorphic function in this half-plane.
Moreover, it is known
\cite{atiyah_part_3}
that the eta function extends meromorphically to the whole $s$-plane
with simple poles.
Formula (10.6) from \cite{jst_part_a} implies
that the eta function does not have a pole at $s=3$ and
that the residue at 
$s=2$ is $4b$,
where $b$ is the coefficient from
(\ref{two-term asymptotic formula for counting function}).


There is an extensive bibliography devoted to the subject of two-term spectral
asymptotics for first order systems. This bibliography spans a period of over
three decades. Unfortunately, all publications prior to \cite{jst_part_a} gave formulae for the
second asymptotic coefficient that were either incorrect or incomplete (i.e. an
algorithm for the calculation of the second asymptotic coefficient rather than
an actual formula). The appropriate bibliographic review is presented in
Section 11 of \cite{jst_part_a}.
The correct explicit formula for the coefficient $b$ is given in Section 1 of \cite{jst_part_a}.

The objective of this paper is to establish the geometric meaning of the coefficient $b$.
The logic behind restricting our analysis to the case when the
manifold is 3-dimensional and $L$ is a $2\times 2$ matrix
differential operator with trace-free principal
symbol and zero subprincipal symbol
is that this is the simplest
eigenvalue
problem for a system of partial differential equations. Hence,
it is ideal for the purpose of establishing the geometric meaning of
the coefficient $b$.

In order to establish the geometric meaning of the coefficient
$b$ we first need to identify the geometric objects encoded
within our eigenvalue problem (\ref{eigenvalue problem}).

\

\textbf{Geometric object 1: the \emph{metric.}}
Observe that the
determinant of the principal symbol is a negative definite quadratic
form in the dual variable (momentum)~$p$,
\begin{equation}
\label{definition of metric}
\det L_\mathrm{prin}(x,p)=-g^{\alpha\beta}p_\alpha p_\beta\,,
\end{equation}
and the coefficients $g^{\alpha\beta}(x)=g^{\beta\alpha}(x)$,
$\alpha,\beta=1,2,3$,
appearing in (\ref{definition of metric}) can be interpreted as the
components of a (contravariant) Riemannian metric.

\

\textbf{Geometric object 2: the \emph{nonvanishing spinor field.}}
The determinant of the principal symbol does
not determine the principal symbol uniquely. In order to identify a
further geometric object encoded within the principal symbol
$L_\mathrm{prin}(x,p)$ we will now start varying this principal symbol,
assuming the metric $g$, defined by formula
(\ref{definition of metric}), to be fixed (prescribed).

Let us fix a reference principal symbol $\mathring L_\mathrm{prin}(x,p)$ corresponding to the
prescribed metric $g$ and look at all principal symbols $L_\mathrm{prin}(x,p)$ which
correspond to the same prescribed metric $g$ and are sufficiently close
to the reference principal symbol.
Restricting our analysis to principal symbols which are close to the reference principal symbol
allows us to avoid dealing with certain topological issues;
this restriction will be dropped in Section ~\ref{Spin structure}.
It turns out, see Section~\ref{Spinor representation of the principal symbol},
that the principal symbols $L_\mathrm{prin}(x,p)$ and $\mathring L_\mathrm{prin}(x,p)$
are related as
\begin{equation}
\label{principal symbol via reference principal symbol}
L_\mathrm{prin}(x,p)=R(x)\,\mathring L_\mathrm{prin}(x,p)\,R^*(x)\,,
\end{equation}
where
\begin{equation}
\label{matrix appearing in unitary transformation of operator}
R:M\to\mathrm{SU}(2)
\end{equation}
is a unique infinitely smooth special unitary matrix-function which is close to the identity matrix.
Thus, special unitary matrix-functions $R(x)$ provide a convenient
parametrization of principal symbols with prescribed metric $g$.

Let $\mathring L$ be the differential operator with principal symbol
$\mathring L_\mathrm{prin}(x,p)$
and zero subprincipal symbol.
It is important to emphasize that for the operators $L$ and $\mathring L$ themselves,
as opposed to their principal symbols, we have, in general, the inequality
\begin{equation}
\label{inequality}
L\ne R\mathring LR^*
\end{equation}
because according to formula (9.3) from \cite{jst_part_a} the
operator $R\mathring LR^*$ has nontrivial subprincipal symbol
$\,\frac i2
\bigl(R_{x^\alpha}(\mathring L_\mathrm{prin})_{p_\alpha}R^*
-
R(\mathring L_\mathrm{prin})_{p_\alpha}R^*_{x^\alpha}\bigr)\,$.
Hence, the transformation of operators $\mathring L\mapsto L$
specified by formula
(\ref{principal symbol via reference principal symbol})
and the conditions that the subprincipal symbols
of $L$ and $\mathring L$ are zero
does, in general,
change the spectrum.

The choice of reference principal symbol $\mathring L_\mathrm{prin}(x,p)$
in our construction is
arbitrary, as long as this principal symbol corresponds to the prescribed metric
$g$, i.e.~as long as we have
$\det\mathring L_\mathrm{prin}(x,p)=-g^{\alpha\beta}(x)\,p_\alpha p_\beta$
for all $(x,p)\in T^*M$.
It is natural to ask the question: what happens if we choose a different reference
principal symbol $\mathring L_\mathrm{prin}(x,p)$?
The freedom in choosing the reference principal symbol $\mathring L_\mathrm{prin}(x,p)$
is a gauge degree of freedom in our construction and our results
are invariant under changes  of the reference principal symbol.
This issue will be addressed in Section~\ref{SU(2) invariance}.

In order to work effectively with special unitary matrices
we need to choose coordinates on the 3-dimensional Lie group $\mathrm{SU}(2)$.
It is convenient to describe a $2\times2$ special unitary matrix by means of a
spinor $\xi$, i.e.~a pair of complex numbers $\xi^a$, $a=1,2$.
The relationship between
a matrix $R\in\mathrm{SU}(2)$ and a nonzero spinor $\xi$ is given by the formula
\begin{equation}
\label{SU(2) matrix expressed via spinor}
R=
\frac1{\|\xi\|}
\begin{pmatrix}
\overline{\xi^1}&\overline{\xi^2}\\
-\xi^2&\xi^1
\end{pmatrix},
\end{equation}
where the overline stands for complex conjugation
and $\|\xi\|:=\sqrt{|\xi^1|^2+|\xi^2|^2}\,$.

Formula (\ref{SU(2) matrix expressed via spinor})
establishes a one-to-one correspondence between
$\mathrm{SU}(2)$ matrices and nonzero spinors,
modulo a rescaling of the spinor by an arbitrary positive real factor.
We choose to specify the scaling of our spinor field $\xi(x)$ in accordance with
\begin{equation}
\label{normalisation of spinor}
\|\xi(x)\|=w(x).
\end{equation}

\begin{remark}
In \cite{jst_part_b} we chose to work with a teleparallel
connection (metric compatible affine connection with zero curvature)
rather than with a spinor field. These are closely
related objects:
locally a teleparallel connection is equivalent to a normalized
(\,$\|\xi(x)\|=1$) spinor field modulo rigid rotations
(\ref{Invariance under rigid rotations equation 3}) of the latter.
\end{remark}

\textbf{Geometric object 3: the \emph{topological charge.}}
It is known, see Section~3 in \cite{jst_part_b}, that the
existence of a principal symbol implies that our manifold $M$ is
parallelizable. Parallelizability, in turn, implies
orientability. Having chosen a particular orientation, we allow
only changes of local coordinates $x^\alpha$, $\alpha=1,2,3$,
which preserve orientation.

We define the topological charge as
\begin{equation}
\label{definition of relative orientation}
\mathbf{c}:=-\frac i2\sqrt{\det g_{\alpha\beta}}\,\operatorname{tr}
\bigl((L_\mathrm{prin})_{p_1}(L_\mathrm{prin})_{p_2}(L_\mathrm{prin})_{p_3}\bigr),
\end{equation}
with the subscripts $p_\alpha$
indicating partial derivatives.
As explained in Section 3 of \cite{jst_part_b},
the number $\mathbf{c}$ defined by formula
(\ref{definition of relative orientation})
can take only two values, $+1$ or $-1$,
and describes the orientation of the principal symbol
relative to the chosen orientation of local coordinates.

Formula (\ref{definition of relative orientation}) defines the topological charge in a purely
analytic fashion. However, later we will give an equivalent definition which is more geometrical,
see formula (\ref{definition of relative orientation more natural}). The frame
$e_j{}^\alpha$ appearing in formula (\ref{definition of relative orientation more natural})
is related to the metric as $g^{\alpha\beta}=\delta^{jk}e_j{}^\alpha e_k{}^\beta$,
so it can be interpreted as the square root of the contravariant metric tensor.
Hence, the topological charge can be loosely described as the sign of the determinant
of the square root of the metric tensor.

\

We have identified three geometric objects encoded within the
eigenvalue problem (\ref{eigenvalue problem})
--- metric, nonvanishing spinor field and topological charge ---
defined in accordance with formulae
(\ref{definition of metric})--(\ref{definition of relative orientation}).
Consequently, one would expect the coefficients $a$ and $b$ from
formula (\ref{two-term asymptotic formula for counting function})
to be expressed via these three geometric objects.
This assertion is confirmed by the following theorem which is
the main result of our paper.

\begin{theorem}
\label{main theorem}
The coefficients in the two-term asymptotics
(\ref{two-term asymptotic formula for counting function})
are given by the formulae
\begin{equation}
\label{formula for a}
a=\frac1{6\pi^2}\int_M
\|\xi\|^3\,\sqrt{\det g_{\alpha\beta}}\ dx\,,
\end{equation}
\begin{equation}
\label{formula for b}
b=\frac{S(\xi)}{2\pi^2}\,,
\end{equation}
where $S(\xi)$ is the massless Dirac action (\ref{definition of massless Dirac action})
with Pauli matrices
\begin{equation}
\label{definition of Pauli matrices sigma}
\sigma^\alpha:=(\mathring L_\mathrm{prin})_{p_\alpha},
\qquad \alpha=1,2,3.
\end{equation}
\end{theorem}

Theorem \ref{main theorem} warrants the following remarks.

Firstly, recall that the $\mathring L$ appearing in Theorem \ref{main theorem}
is our reference operator which we need to describe all
possible operators $L$ with given metric $g$.
What happens if we take $L=\mathring L$? In this case
formula (\ref{formula for b}) holds with spinor field
$\xi^1(x)=w(x)$, $\xi^2(x)=0$. This, on its own, is a
nontrivial result.

Secondly, the topological charge $\mathbf{c}$ does not appear
explicitly in Theorem~\ref{main theorem}. Nevertheless,
it is implicitly present in our Pauli matrices
(\ref{definition of Pauli matrices sigma}).
Indeed, formula
(\ref{principal symbol via reference principal symbol})
implies that the integer quantity
\[
-\frac i2\sqrt{\det g_{\alpha\beta}}\,\operatorname{tr}
\bigl(
(\mathring L_\mathrm{prin})_{p_1}
(\mathring L_\mathrm{prin})_{p_2}
(\mathring L_\mathrm{prin})_{p_3}
\bigr)
\]
has the same value as (\ref{definition of relative orientation}).

Thirdly, it is tempting to apply Theorem~\ref{main theorem} in
the case when the operator $L$ is itself a massless Dirac
operator. This cannot be done because a massless Dirac
operator acts on spinors rather than on pairs of
half-densities. This impediment can be overcome by switching to
a massless Dirac operator on half-densities, see formula (A.19)
in \cite{jst_part_b}. However, we cannot take $L$ to be a
massless Dirac operator on half-densities either because, according to
Lemma 6.1 from \cite{jst_part_b}, the latter has a nontrivial
subprincipal symbol.
Furthermore, it is known \cite{bismut_and_freed,jst_part_b}
that for the massless Dirac
operator the coefficient $b$ is zero.

Finally, Theorem~\ref{main theorem} provides a fresh perspective
on the history of the subject of two-term spectral asymptotics
for first order systems, see Section 11 of~\cite{jst_part_a} for
details. Namely, Theorem~\ref{main theorem} shows that even in
the simplest case the second asymptotic coefficient for a first
order system has a highly nontrivial geometric meaning. At a
formal level, the application of microlocal techniques does not
require the use of advanced differential geometric concepts.
However, the calculations involved are so complicated that it is
hard to avoid mistakes without an understanding of the
differential geometric content of the spectral problem.

It is also worth noting that we use the term `Pauli matrices' in a more
general sense than in traditional quantum mechanics.
The traditional definition is the one from formula (\ref{Pauli matrices s})
and it corresponds to flat space,
whereas our definition is adapted to curved space.
For us Pauli matrices $\sigma^\alpha$ are trace-free Hermitian $2\times2$ matrices satisfying the
identity (\ref{defining relation for Pauli matrices}).
It might have been be more appropriate to call our matrices
$\sigma^\alpha$, $\alpha=1,2,3$, Pauli matrices of Riemannian metric $g^{\alpha\beta}$,
but, as this expression is too long, we call them simply Pauli matrices.

The paper is organized as follows.
In Section~\ref{Spinor representation of the principal symbol}
we explain the origins of formula
(\ref{principal symbol via reference principal symbol})
and
in Section~\ref{Proof of Theorem} we give
the proof of Theorem~\ref{main theorem}.
In Section~\ref{Spin structure}
we introduce the concept of spin structure
which allows us to drop the restriction that
our principal symbol $L_\mathrm{prin}(x,p)$
is sufficiently close to the reference principal symbol
$\mathring L_\mathrm{prin}(x,p)$.
And in Sections~\ref{Conformal invariance}--\ref{Invariance under rigid rotations}
we show that our formula
(\ref{formula for b})
is invariant under the action of certain gauge transformations.

\section{Spinor representation of the principal symbol}
\label{Spinor representation of the principal symbol}

Let $L_\mathrm{prin}(x,p)$ and $\mathring L_\mathrm{prin}(x,p)$
be a pair of trace-free Hermitian $2\times2$ principal symbols
and let $g$ be a prescribed Riemannian metric.
Both $L_\mathrm{prin}(x,p)$ and $\mathring L_\mathrm{prin}(x,p)$
are assumed to be linear in $p$:
\begin{eqnarray}
\label{principal symbol of A is linear in xi}
L_\mathrm{prin}(x,p)&=&L_\mathrm{prin}^{(\alpha)}(x)\,p_\alpha\,,
\\
\label{principal symbol of B is linear in xi}
\mathring L_\mathrm{prin}(x,p)&=&\mathring L_\mathrm{prin}^{(\alpha)}(x)\,p_\alpha\,,
\end{eqnarray}
where $L_\mathrm{prin}^{(\alpha)}(x)$ and $\mathring L_\mathrm{prin}^{(\alpha)}(x)$,
$\alpha=1,2,3$, are some trace-free Hermitian
$2\times2$ matrix-functions.
The assumption that our principal symbols
$L_\mathrm{prin}(x,p)$ and $\mathring L_\mathrm{prin}(x,p)$ are linear in $p$
means, of course, that we are dealing with differential operators as opposed
to pseudodifferential operators.
The principal symbols
$L_\mathrm{prin}(x,p)$ and $\mathring L_\mathrm{prin}(x,p)$ are assumed to satisfy
\begin{equation}
\label{definition of metric with A and B}
\det L_\mathrm{prin}(x,p)=\det\mathring L_\mathrm{prin}(x,p)
=-g^{\alpha\beta}(x)\,p_\alpha p_\beta
\end{equation}
for all $(x,p)\in T^*M$,
and are also assumed to be sufficiently close in terms of the $C^\infty(M)$ topology
applied to the matrix-functions
$L_\mathrm{prin}^{(\alpha)}(x)$ and $\mathring L_\mathrm{prin}^{(\alpha)}(x)$,
$\alpha=1,2,3$.

Our task in this section is to show that there exists a unique infinitely smooth
special unitary matrix-function
(\ref{matrix appearing in unitary transformation of operator})
which is close to the identity matrix and which relates the principal
symbols $L_\mathrm{prin}(x,p)$ and $\mathring L_\mathrm{prin}(x,p)$
in accordance with formula (\ref{principal symbol via reference principal symbol}).

We follow the convention of \cite{MR2670535,jst_part_b} in
denoting the elements of the matrices
$L_\mathrm{prin}^{(\alpha)}$ and $\mathring L_\mathrm{prin}^{(\alpha)}$
as
$\bigl(L_\mathrm{prin}^{(\alpha)}\bigr)_{\dot ab}$
and
$\bigl(\mathring L_\mathrm{prin}^{(\alpha)}\bigr)_{\dot ab}$
respectively,
where the dotted index, running through the values $\dot1,\dot2$, enumerates the rows
and the undotted index, running through the values $1,2$, enumerates the columns.

Put
\begin{eqnarray}
\label{frame via principal symbol of A}
e_1{}^\alpha:=\operatorname{Re}\bigl(L_\mathrm{prin}^{(\alpha)}\bigr)_{\dot12}\,,
\quad
e_2{}^\alpha:=-\operatorname{Im}\bigl(L_\mathrm{prin}^{(\alpha)}\bigr)_{\dot12}\,,
\quad
e_3{}^\alpha:=\operatorname{Re}\bigl(L_\mathrm{prin}^{(\alpha)}\bigr)_{\dot11}\,,
\\
\label{frame via principal symbol of B}
\mathring e_1{}^\alpha:=\operatorname{Re}
\bigl(\mathring L_\mathrm{prin}^{(\alpha)}\bigr)_{\dot12}\,,
\quad
\mathring e_2{}^\alpha:=-\operatorname{Im}
\bigl(\mathring L_\mathrm{prin}^{(\alpha)}\bigr)_{\dot12}\,,
\quad
\mathring e_3{}^\alpha:=\operatorname{Re}
\bigl(\mathring L_\mathrm{prin}^{(\alpha)}\bigr)_{\dot11}\,.
\end{eqnarray}
As explained in Section 3 of \cite{jst_part_b},
formula (\ref{frame via principal symbol of A}) defines a \emph{frame} ---
a triple of infinitely smooth real orthonormal vector fields
$e_j(x)$, $j=1,2,3$, on the manifold $M$ ---
and, moreover, the principal symbol $L_\mathrm{prin}(x,p)$
is equivalent to the frame $e_j$ in the sense that the
principal symbol uniquely determines the frame and
the frame uniquely determines the principal symbol.
Similarly, formula (\ref{frame via principal symbol of B}) defines a frame
$\mathring e_j$ which is equivalent to the principal symbol $\mathring L_\mathrm{prin}(x,p)$.

Condition (\ref{definition of metric with A and B})
implies that the frames $e_j$ and $\mathring e_j$ are orthonormal
with respect to the same metric. Hence,
this pair of orthonormal frames is related as
\begin{equation}
\label{orthogonal transformation of frame}
e_j{}^\alpha=O_j{}^k\mathring e_k{}^\alpha,
\end{equation}
where $O(x)$ is a $3\times3$ orthogonal matrix-function
with elements
\[
O_j{}^k=\delta^{kl}g_{\alpha\beta}\,e_j{}^\alpha\,\mathring e_l{}^\beta.
\]
As we assumed the principal symbols
$L_\mathrm{prin}(x,p)$ and $\mathring L_\mathrm{prin}(x,p)$ to be close,
the frames $e_j$ and $\mathring e_j$ are also close. Consequently, the matrix-function
$O(x)$ is close to the identity matrix and, hence, special orthogonal.

It is well-known that the Lie group $\mathrm{SO}(3)$ is locally (in a neighbourhood of the identity)
isomorphic to the Lie group $\mathrm{SU}(2)$.
According to formulae (A.15) and (A.2) from \cite{jst_part_b},
a $3\times3$ special orthogonal matrix $O$
is expressed via a $2\times2$ special unitary matrix $R$ as
\begin{equation}
\label{special orthogonal matrix via special unitary matrix}
O_j{}^k=\frac12\operatorname{tr}(s_jRs^kR^*)\,,
\end{equation}
where
\begin{equation}
\label{Pauli matrices s}
s^1:=
\begin{pmatrix}
0&1\\
1&0
\end{pmatrix}
=s_1
\,,
\quad
s^2:=
\begin{pmatrix}
0&-i\\
i&0
\end{pmatrix}
=s_2
\,,
\quad
s^3:=
\begin{pmatrix}
1&0\\
0&-1
\end{pmatrix}
=s_3
\,.
\end{equation}
Formula (\ref{special orthogonal matrix via special unitary matrix}) tells us that
a $3\times3$ special orthogonal matrix is, effectively,
the square of a $2\times2$ special unitary matrix.
Formula (\ref{special orthogonal matrix via special unitary matrix}) provides a local
diffeomorphism between neighbourhoods of the identity
in $\mathrm{SO}(3)$ and in $\mathrm{SU}(2)$.

A straightforward calculation shows that formulae
(\ref{principal symbol of A is linear in xi}),
(\ref{principal symbol of B is linear in xi})
and
(\ref{frame via principal symbol of A})--(\ref{Pauli matrices s})
imply formula (\ref{principal symbol via reference principal symbol}).

Let us now define Pauli matrices $\sigma^\alpha$ in accordance with formula
(\ref{definition of Pauli matrices sigma}). Of course, we have
\begin{equation}
\label{definition of Pauli matrices sigma alternative}
\sigma^\alpha(x)=\mathring L_\mathrm{prin}^{(\alpha)}(x),
\qquad \alpha=1,2,3,
\end{equation}
where the $\mathring L_\mathrm{prin}^{(\alpha)}$
are the matrix-functions from formula
(\ref{principal symbol of B is linear in xi}).
We could stick with the notation $\mathring L_\mathrm{prin}^{(\alpha)}$
but we choose to switch
to $\sigma^\alpha$ because this is how Pauli matrices are traditionally denoted in the subject.

It is easy to see that formula (\ref{definition of metric with A and B}) implies
\begin{equation}
\label{defining relation for Pauli matrices}
\sigma^\alpha\sigma^\beta+\sigma^\beta\sigma^\alpha=2Ig^{\alpha\beta},
\end{equation}
where $I$ is the $2\times2$ identity matrix.
Formula  (\ref{defining relation for Pauli matrices})
means that our $\sigma^\alpha$
satisfy the defining relation for Pauli matrices.

Formulae
(\ref{orthogonal transformation of frame})--(\ref{Pauli matrices s}),
(\ref{SU(2) matrix expressed via spinor}),
(\ref{frame via principal symbol of B})
and
(\ref{definition of Pauli matrices sigma alternative})
allow us to express the frame $e_j$ via the spinor field $\xi$ and Pauli
matrices $\sigma^\alpha$.
We took great care to choose coordinates on the Lie group $\mathrm{SU}(2)$
(i.e.~structure of the matrix in the RHS of formula (\ref{SU(2) matrix expressed via spinor}))
so that the resulting expressions agree with formulae (B.3), (B.4) and (B.1)
from~\cite{MR2670535}. The only difference is in notation:
the $\vartheta^j$ in Appendix B of \cite{MR2670535} stands for
$\vartheta^j{}_\alpha=\delta^{jk}g_{\alpha\beta}\,e_k{}^\beta$
(compare with formula (\ref{most basic definition of coframe})).

The fact that our construction agrees with that in \cite{MR2670535}
will become important in the next section when we will make use
of a particular formula from \cite{MR2670535}.

\begin{remark}
\label{remark about topological charge}
As explained in Section~3 of \cite{jst_part_b},
the topological charge, initially defined in accordance with formula
(\ref{definition of relative orientation}),
can be equivalently rewritten in terms of frames as
\begin{equation}
\label{definition of relative orientation more natural}
\mathbf{c}=\operatorname{sgn}\det e_j{}^\alpha
=\operatorname{sgn}\det\mathring e_j{}^\alpha.
\end{equation}
The paper \cite{MR2670535} was written under the assumption that
\begin{equation}
\label{topological invariant equals plus one}
\mathbf{c}=+1\,,
\end{equation}
see formula (A.1) in \cite{MR2670535}.
This means that care is required when using the results of \cite{MR2670535}.
Namely, in the next section we will first prove Theorem \ref{main theorem}
for the case (\ref{topological invariant equals plus one}) and then provide a separate
argument explaining why formula (\ref{formula for b}) remains true in the case
\begin{equation}
\label{topological invariant equals minus one}
\mathbf{c}=-1\,.
\end{equation}
\end{remark}

\section{Proof of Theorem \ref{main theorem}}
\label{Proof of Theorem}

We prove Theorem \ref{main theorem} by examining the equivalent spectral
problem (\ref{eigenvalue problem without weight}). Note that transition from
(\ref{eigenvalue problem})
to
(\ref{eigenvalue problem without weight})
is a special case of the gauge transformation
(\ref{gauge transformation of eigenvalue problem})
with $\varphi=\ln w$.
As explained in the beginning of Section~\ref{Conformal invariance},
the transformation
(\ref{gauge transformation of eigenvalue problem})
preserves the structure of our eigenvalue problem:
the principal symbol of the operator $w^{-1/2}Lw^{-1/2}$
is trace-free and its subprincipal symbol is zero.

We now apply Theorem 1.1 from \cite{jst_part_b} to the
eigenvalue problem (\ref{eigenvalue problem without weight}).

Our formula (\ref{formula for a}) is an immediate consequence of
formula (1.18) from \cite{jst_part_b}
and our formulae (\ref{definition of metric}) and (\ref{normalisation of spinor}).
Here, of course, we use the fact that we are working in dimension three.

The proof of formula (\ref{formula for b}) is more delicate so we
initially consider the case
\begin{equation}
\label{case of trivial weight function}
w(x)=1,\qquad\forall x\in M.
\end{equation}
In this case, according to formulae (1.19) and (8.1) from \cite{jst_part_b}, we have
\begin{equation}
\label{formula for b when w is 1 meaningful}
b=\frac{3\mathbf{c}}{8\pi^2}\int_M
*T^\mathrm{ax}
\,\sqrt{\det
g_{\alpha\beta}}\ dx\,,
\end{equation}
where
\begin{multline}
\label{explicit formula for the trace of torsion with a star}
*\!T^\mathrm{ax}
=\frac{\delta_{kl}}3\,\sqrt{\det g^{\alpha\beta}}
\ \bigl[
e^k{}_{1}
\,\partial e^l{}_{3}/\partial x^2
+
e^k{}_{2}
\,\partial e^l{}_{1}/\partial x^3
+
e^k{}_{3}
\,\partial e^l{}_{2}/\partial x^1
\\
-
e^k{}_{1}
\,\partial e^l{}_{2}/\partial x^3
-
e^k{}_{2}
\,\partial e^l{}_{3}/\partial x^1
-
e^k{}_{3}
\,\partial e^l{}_{1}/\partial x^2
\bigr],
\end{multline}
\begin{equation}
\label{most basic definition of coframe}
e^j{}_\alpha=\delta^{jk}g_{\alpha\beta}\,e_k{}^\beta.
\end{equation}
The real scalar field $\,*T^\mathrm{ax}(x)\,$ has the geometric
meaning of the Hodge dual of axial torsion of the teleparallel connection,
see \cite{jst_part_b} for details.

Let us now drop the assumption
(\ref{case of trivial weight function}).

The introduction of a weight function is equivalent to a scaling of the principal
symbol
$L_\mathrm{prin}(x,p)\mapsto(w(x))^{-1}L_\mathrm{prin}(x,p)$,
which, in view of formulae
(\ref{frame via principal symbol of A})
and
(\ref{definition of metric}),
leads to a scaling of the frame
\begin{equation}
\label{scaling of frame in proof of theorem}
e_j\mapsto w^{-1}e_j
\end{equation}
and corresponding scaling of the metric
\begin{equation}
\label{scaling of metric in proof of theorem}
g^{\alpha\beta}\mapsto w^{-2}g^{\alpha\beta}.
\end{equation}
Substituting
(\ref{scaling of frame in proof of theorem})
and
(\ref{scaling of metric in proof of theorem})
into
(\ref{most basic definition of coframe})
and
(\ref{explicit formula for the trace of torsion with a star})
we see that the integrand in formula
(\ref{formula for b when w is 1 meaningful})
scales as
\begin{equation}
\label{rescaling of integrand}
*\!T^\mathrm{ax}\,\sqrt{\det g_{\alpha\beta}}
\ \mapsto
\ w^2*\!T^\mathrm{ax}\,\sqrt{\det g_{\alpha\beta}}\ .
\end{equation}
Here the remarkable fact is that we do not get derivatives of the weight
function because these cancel out due to the antisymmetric
structure of the RHS of formula
(\ref{explicit formula for the trace of torsion with a star}).
In other words, axial torsion,
defined by formulae (1.20) and (3.12) from~\cite{jst_part_b},
has the remarkable property that it
scales in a covariant manner under scaling of the frame.
Note that the full torsion tensor,
defined by formula (3.12) from~\cite{jst_part_b},
does not possess such a covariance property.

Formula (\ref{rescaling of integrand}) tells us that in order to accommodate an arbitrary
weight function $w(x)$ we need to multiply the integrand in formula
(\ref{formula for b when w is 1 meaningful}) by $(w(x))^2$, which gives us
\begin{equation}
\label{formula for b when w is 1 meaningful rescaled}
b=\frac{3\mathbf{c}}{8\pi^2}\int_M
w^2*\!T^\mathrm{ax}
\,\sqrt{\det
g_{\alpha\beta}}\ dx\,.
\end{equation}
Let us emphasize that the metric and torsion appearing in formula
(\ref{formula for b when w is 1 meaningful rescaled})
are the original, unscaled metric and torsion determined by the
original, unscaled principal symbol $L_\mathrm{prin}(x,p)$.
The scaling argument has been incorporated into the factor $(w(x))^2$.

We now need to express the integrand in
(\ref{formula for b when w is 1 meaningful rescaled})
in terms of the spinor field $\xi$.

We already have an expression for the weight function in terms of the spinor field,
see formula (\ref{normalisation of spinor}). So we only need to express the Hodge
dual of axial torsion in terms of the spinor field.
Formulae
(\ref{frame via principal symbol of A}),
(\ref{principal symbol of A is linear in xi}),
(\ref{principal symbol of B is linear in xi}),
(\ref{definition of Pauli matrices sigma alternative}),
(\ref{principal symbol via reference principal symbol})
and
(\ref{SU(2) matrix expressed via spinor})
allow us to express the frame $e_j$
via the spinor field $\xi$ and Pauli
matrices $\sigma^\alpha$.
Hence one needs to combine all these formulae to
get explicit expressions for the vector fields $e_j$, $j=1,2,3$, and
substitute these into
(\ref{most basic definition of coframe})
and
(\ref{explicit formula for the trace of torsion with a star}).
This is a massive calculation. Fortunately,
for the case
(\ref{topological invariant equals plus one})
this calculation was carried out
in Appendix B of \cite{MR2670535}: formula
(B.5) from \cite{MR2670535} reads
\begin{equation}
\label{formula from MR2670535}
*\!T^\mathrm{ax}=\frac{4\operatorname{Re}(\xi^*W\xi)}{3\|\xi\|^2}\,,
\end{equation}
where $W$ is the massless Dirac operator (\ref{definition of Weyl operator}).

Formulae
(\ref{formula for b when w is 1 meaningful rescaled}),
(\ref{normalisation of spinor}),
(\ref{formula from MR2670535})
and
(\ref{definition of massless Dirac action})
imply formula
(\ref{formula for b}).
This completes the proof of Theorem \ref{main theorem}
for the case (\ref{topological invariant equals plus one}).

In order to prove formula (\ref{formula for b}) for the case
(\ref{topological invariant equals minus one}),
we invert coordinates ($x^\alpha\mapsto-x^\alpha$),
which changes the sign of topological charge
and allows us to use formula (\ref{formula for b}).
We then invert coordinates again and use the facts that
\begin{itemize}
\item
the integrand of the massless Dirac action
(\ref{definition of massless Dirac action})
is invariant under inversion of coordinates
and
\item
our spinor field $\xi$ defined by formulae
(\ref{principal symbol via reference principal symbol})--(\ref{normalisation of spinor})
is an anholonomic object,
i.e.~it does not depend on the choice of coordinates.
\end{itemize}

\section{Spin structure}
\label{Spin structure}

In stating our results in Section~\ref{Main result}
we assumed the principal symbols
$L_\mathrm{prin}(x,p)$ and $\mathring L_\mathrm{prin}(x,p)$
to be sufficiently close. This was done in order to ensure that
equation (\ref{principal symbol via reference principal symbol})
could be resolved with respect to the special unitary matrix-function
$R(x)$. The restriction of closeness of principal symbols can be overcome
by means of the introduction of the concept of spin structure.

\begin{definition}
\label{definition of spin structure}
We say that the principal symbols
$L_\mathrm{prin}(x,p)$ and $\mathring L_\mathrm{prin}(x,p)$
have the same
spin structure if there exists an infinitely smooth special unitary matrix-function
(\ref{matrix appearing in unitary transformation of operator})
such that we have (\ref{principal symbol via reference principal symbol}).
\end{definition}

\begin{remark}
Principal symbols  with the same
spin structure form an equivalence class.
\end{remark}

The closeness of the principal symbols
$L_\mathrm{prin}(x,p)$ and $\mathring L_\mathrm{prin}(x,p)$
was never
used in the proof of Theorem~\ref{main theorem}. All that is needed for
Theorem~\ref{main theorem} to be true is for
the principal symbols
$L_\mathrm{prin}(x,p)$ and $\mathring L_\mathrm{prin}(x,p)$
to have the same spin structure,
i.e.~belong to the same equivalence class.

Hence, it would have been more logical
to identify the spin structure as a separate geometric object from the very
start, in Section~\ref{Main result}, and avoid arguments relying on the
closeness of the principal symbols. We chose not to proceed along this route
in order to make the exposition in Section~\ref{Main result} as simple and clear as possible.

The only difference between the ``local'' setting
(the principal symbols
$L_\mathrm{prin}(x,p)$ and $\mathring L_\mathrm{prin}(x,p)$
are assumed to be close)
and the ``global'' setting
(the principal symbols
$L_\mathrm{prin}(x,p)$ and $\mathring L_\mathrm{prin}(x,p)$
are assumed to
have the same spin structure) is that we can no longer
claim that the special unitary matrix-function $R(x)$ appearing in formula
(\ref{principal symbol via reference principal symbol}) is unique.
In the ``local'' setting uniqueness was achieved by requiring $R(x)$
to be close to the identity matrix, whereas in the ``global'' setting
$R(x)$ is defined modulo sign
(not surprising as $\mathrm{SU}(2)$ is the double cover of $\mathrm{SO}(3)$).
This  sign indeterminacy does not affect formula
(\ref{formula for b})
because the massless Dirac action is quadratic in the spinor field.

The number of different spin structures (i.e.~number of equivalence classes
of principal symbols) depends on the topology of the manifold. Say,
the torus $\mathbb{T}^3$ admits eight different spin structures,
whereas the sphere $\mathbb{S}^3$ admits a unique spin structure.
See Appendices A and B in \cite{jst_part_b} and further bibliographic
references therein for more details.

It may seem that our
Definition~\ref{definition of spin structure}
is different from the definition of spin structure
in differential geometric literature. Indeed, differential geometers do
not operate with concepts such as the principal symbol, using
frames instead. However, it has been shown in Section 3 of \cite{jst_part_b}
that a principal symbol is equivalent to a frame, so
our ``microlocal'' definition of spin structure
can be easily recast in terms of frames, bringing it in agreement
with the traditional differential geometric one.

Here it is important to emphasize that we do not claim to have
redefined the notion of spin structure for the most general
case. We work in the very specific setting of dimension
three.

\section{Conformal invariance}
\label{Conformal invariance}

Let us transform the differential operator $L$ and weight $w(x)$ as
\begin{equation}
\label{gauge transformation of eigenvalue problem}
L\mapsto e^{-\varphi/2}Le^{-\varphi/2},
\qquad
w\mapsto e^{-\varphi}w,
\end{equation}
where $\varphi:M\to\mathbb{R}$ is an arbitrary infinitely smooth
real-valued scalar function.
The transformation (\ref{gauge transformation of eigenvalue problem})
does not change the spectrum of our eigenvalue problem
(\ref{eigenvalue problem}) and, moreover, preserves its structure:
the principal symbol remains trace-free and the subprincipal symbol
remains zero. The fact that the
subprincipal symbol remains zero is established by using formula
(9.3) from \cite{jst_part_a}. Here, of course, it is important that
we are conjugating the operator by a real-valued scalar function
$e^{-\varphi/2}$
rather than a complex-valued matrix-function $R$.

In this section we examine how the gauge transformation
(\ref{gauge transformation of eigenvalue problem})
works its way into scalings of the metric, Pauli matrices and spinor field.

Formulae (\ref{definition of metric}) and (\ref{gauge transformation of eigenvalue problem})
imply that the metric transforms as
\begin{equation}
\label{scaling of metric}
g^{\alpha\beta}\mapsto e^{-2\varphi}g^{\alpha\beta},
\end{equation}
which means that we are looking at a conformal scaling of the metric.

We scale the reference principal symbol $\mathring L_\mathrm{prin}(x,p)$ the same way
as the principal symbol $L_\mathrm{prin}(x,p)$,
i.e.~as
\begin{equation}
\label{scaling of reference principal symbol}
\mathring L_\mathrm{prin}\mapsto e^{-\varphi}\mathring L_\mathrm{prin}\,,
\end{equation}
because this way we maintain the condition
(\ref{definition of metric with A and B}).
Formulae (\ref{definition of Pauli matrices sigma})
and
(\ref{scaling of reference principal symbol})
imply that the Pauli matrices scale as
\begin{equation}
\label{scaling of Pauli matrices}
\sigma^\alpha\mapsto e^{-\varphi}\sigma^\alpha.
\end{equation}

Of course, the scaling of Pauli matrices
(\ref{scaling of Pauli matrices})
agrees with the scaling of the metric
(\ref{scaling of metric})
in the sense that the scaled
Pauli matrices and metric satisfy the identity
(\ref{defining relation for Pauli matrices}).

Formulae
(\ref{normalisation of spinor})
and
(\ref{gauge transformation of eigenvalue problem})
imply that the spinor field scales as
\begin{equation}
\label{scaling of spinor}
\xi\mapsto e^{-\varphi}\xi.
\end{equation}

Let us now examine what happens to the massless Dirac action
(\ref{definition of massless Dirac action}) under the transformations
(\ref{scaling of metric}),
(\ref{scaling of Pauli matrices})
and
(\ref{scaling of spinor}).

We first look at the expression $W\xi$.
Examination of formulae
(\ref{definition of Weyl operator})
and
(\ref{definition of Christoffel symbols})
show that the expression $W\xi$ transforms as
\begin{equation}
\label{transformation of W zeta}
W\xi\mapsto e^{-2\varphi}W\xi.
\end{equation}
We see that the expression $W\xi$ scales in a covariant way,
with ``covariant'' meaning that the derivatives of $\varphi$ do not appear
in the RHS of (\ref{transformation of W zeta}). Of course, the covariance
of the massless Dirac operator under conformal scalings of the metric
is a know differential geometric fact: see Theorem 4.3 in \cite{solovej}.

Formulae (\ref{scaling of spinor}) and (\ref{transformation of W zeta}) imply
\begin{equation}
\label{transformation of massless dirac lagrangian without density}
\operatorname{Re}(\xi^*W\xi)\mapsto e^{-3\varphi}\operatorname{Re}(\xi^*W\xi).
\end{equation}
Formula (\ref{scaling of metric}) implies
$g_{\alpha\beta}\mapsto e^{2\varphi}g_{\alpha\beta}$ and,
as we are working in dimension 3, this, in turn, implies that the Riemannian density scales as
\begin{equation}
\label{scaling of density}
\sqrt{\det g_{\alpha\beta}}\mapsto e^{3\varphi}\,\sqrt{\det g_{\alpha\beta}}\,.
\end{equation}
Substituting formulae
(\ref{transformation of massless dirac lagrangian without density})
and
(\ref{scaling of density})
into formula
(\ref{definition of massless Dirac action})
we see that our massless Dirac action is invariant
under the transformations
(\ref{scaling of metric})
(\ref{scaling of Pauli matrices})
and
(\ref{scaling of spinor}).
This is, of course, in agreement with Theorem \ref{main theorem}:
as the gauge transformation (\ref{gauge transformation of eigenvalue problem})
does not change the spectrum of our eigenvalue problem (\ref{eigenvalue problem}),
it does not change the second asymptotic coefficient
(\ref{formula for b})
of the counting function.

\section{$\mathrm{SU}(2)$ invariance}
\label{SU(2) invariance}

Let us transform the reference principal symbol $\mathring L_\mathrm{prin}(x,p)$ as
\begin{equation}
\label{SU(2) invariance equation 1}
\mathring L_\mathrm{prin}\mapsto Q\mathring L_\mathrm{prin}Q^*,
\end{equation}
where $Q:M\to\mathrm{SU}(2)$ is an arbitrary infinitely smooth special unitary
matrix-function.
Formulae
(\ref{definition of Pauli matrices sigma})
and
(\ref{SU(2) invariance equation 1})
imply
\begin{equation}
\label{SU(2) invariance equation 2}
\sigma^\alpha\mapsto Q\sigma^\alpha Q^*.
\end{equation}
Also, formulae
(\ref{principal symbol via reference principal symbol})
and
(\ref{SU(2) invariance equation 1})
imply $R\mapsto RQ^*$, which can be equivalently rewritten as
\begin{equation}
\label{SU(2) invariance equation 3}
R^*\mapsto QR^*.
\end{equation}
Examining the structure of the matrix $R$, see formula
(\ref{SU(2) matrix expressed via spinor}),
we conclude that formula
(\ref{SU(2) invariance equation 3})
is equivalent to the linear transformation of the spinor field
\begin{equation}
\label{SU(2) invariance equation 4}
\xi\mapsto Q\xi.
\end{equation}
Formulae
(\ref{SU(2) invariance equation 2}),
(\ref{SU(2) invariance equation 4})
and Property 4 from Appendix A of \cite{jst_part_b}
tell us that our massless Dirac action  is invariant
under the transformation
(\ref{SU(2) invariance equation 1}).
This is, of course, in agreement with Theorem \ref{main theorem}:
the choice of reference principal symbol
does not affect the spectrum of our eigenvalue problem (\ref{eigenvalue problem}),
hence, it does not affect the second asymptotic coefficient
(\ref{formula for b})
of the counting function.

\section{Invariance under rigid rotations}
\label{Invariance under rigid rotations}

Let us transform the differential operator $L$ as
\begin{equation}
\label{Invariance under rigid rotations equation 1}
L\mapsto\mathbf{Q}L\mathbf{Q}^*,
\end{equation}
where
$
\mathbf{Q}=
\begin{pmatrix}
\mathbf{Q}_{11}&\mathbf{Q}_{12}
\\
\mathbf{Q}_{21}&\mathbf{Q}_{22}
\end{pmatrix}
$
is a \textbf{constant} special unitary matrix. The transformation
(\ref{Invariance under rigid rotations equation 1})
does not change the spectrum of our eigenvalue problem
(\ref{eigenvalue problem}) and, moreover, preserves its structure:
the principal symbol remains trace-free and the subprincipal symbol
remains zero. We refer to the transformation
(\ref{Invariance under rigid rotations equation 1})
as a \emph{rigid rotation} because
it describes a rotation of the frame (\ref{frame via principal symbol of A}),
with this rotation being the same at all points of the manifold $M$.

The transformation (\ref{Invariance under rigid rotations equation 1})
is equivalent to the following transformation of the
special unitary matrix-function $R(x)$ appearing in formula
(\ref{principal symbol via reference principal symbol}):
\begin{equation}
\label{Invariance under rigid rotations equation 2}
R\mapsto\mathbf{Q}R.
\end{equation}
Formulae
(\ref{SU(2) matrix expressed via spinor})
and
(\ref{normalisation of spinor})
give us a one-to-one correspondence
between special unitary matrix-functions and weight functions on the one hand
and nonvanishing spinor fields on the other.
In terms of the spinor field
the transformation (\ref{Invariance under rigid rotations equation 2})
reads
\begin{equation}
\label{Invariance under rigid rotations equation 3}
\begin{pmatrix}
\xi^1
\\
\xi^2
\end{pmatrix}
\mapsto
\begin{pmatrix}
\mathbf{Q}_{21}\overline{\xi^2}+\mathbf{Q}_{22}\xi^1
\\
-\mathbf{Q}_{21}\overline{\xi^1}+\mathbf{Q}_{22}\xi^2
\end{pmatrix}.
\end{equation}
Note that, unlike (\ref{SU(2) invariance equation 3}),
this transformation is \textbf{not} linear because of the complex conjugation.
The transformation
(\ref{Invariance under rigid rotations equation 3}) can be written
as a sum of linear and antilinear transformations:
\begin{equation}
\label{Invariance under rigid rotations equation 4}
\xi\mapsto
\mathbf{Q}_{22}\xi-\mathbf{Q}_{21}\mathrm{C}(\xi)
\end{equation}
where
$
\begin{pmatrix}
\xi^1
\\
\xi^2
\end{pmatrix}
\mapsto
\mathrm{C}(\xi):=
\begin{pmatrix}
-\overline{\xi^2}
\\
\overline{\xi^1}
\end{pmatrix}
$
is the charge conjugation operator, see formula (A.9)
in \cite{jst_part_b}.

Let us show, by performing explicit calculations, that
the massless Dirac action
(\ref{definition of massless Dirac action})
is invariant under the transformation
(\ref{Invariance under rigid rotations equation 4}).
Using the fact that the massless Dirac operator commutes with the
charge conjugation operator,
see Property 3 in Appendix A of \cite{jst_part_b},
we get
\begin{multline*}
(\mathbf{Q}_{22}\xi-\mathbf{Q}_{21}\mathrm{C}(\xi))^*
\,W\,
(\mathbf{Q}_{22}\xi-\mathbf{Q}_{21}\mathrm{C}(\xi))
\\
=|\mathbf{Q}_{22}|^2\,
\xi^*W\xi
+|\mathbf{Q}_{21}|^2\,
(\mathrm{C}(\xi))^*\mathrm{C}(W\xi)
-\overline{\mathbf{Q}_{22}}\mathbf{Q}_{21}
\xi^*\mathrm{C}(W\xi)
-\mathbf{Q}_{22}\overline{\mathbf{Q}_{21}}
(\mathrm{C}(\xi))^*W\xi
\\
=|\mathbf{Q}_{22}|^2\,
\xi^*W\xi
+|\mathbf{Q}_{21}|^2\,
\overline{\xi^*W\xi}
+\bigl\{
\overline{\mathbf{Q}_{22}}\mathbf{Q}_{21}
(\mathrm{C}(\xi))^T\,\overline{W\xi}\,
-\mathbf{Q}_{22}\overline{\mathbf{Q}_{21}}
(\mathrm{C}(\xi))^*W\xi
\bigr\}.
\end{multline*}
But the expression in the curly brackets is purely imaginary, so
\begin{multline*}
\operatorname{Re}
\bigl[
(\mathbf{Q}_{22}\xi-\mathbf{Q}_{21}\mathrm{C}(\xi))^*
\,W\,
(\mathbf{Q}_{22}\xi-\mathbf{Q}_{21}\mathrm{C}(\xi))
\bigr]
\\
=|\mathbf{Q}_{22}|^2
\operatorname{Re}
\bigl[
\xi^*W\xi
\bigr]
+|\mathbf{Q}_{21}|^2
\operatorname{Re}
\bigl[
\,\overline{\xi^*W\xi}\,
\bigr]
=\bigl(
|\mathbf{Q}_{22}|^2
+
|\mathbf{Q}_{21}|^2
\bigr)
\operatorname{Re}
\bigl[
\xi^*W\xi
\bigr]
=\operatorname{Re}
\bigl[
\xi^*W\xi
\bigr].
\end{multline*}

\appendix

\section{Invariant analytic description of a first order differential operator}
\label{Invariant analytic description of a first order differential operator}

Let $L$ be a formally self-adjoint first order linear differential
operator acting on $m$-columns
$v=\begin{pmatrix}v_1&\ldots&v_m\end{pmatrix}^T$
of complex-valued half-densities
over a connected $n$-dimensional manifold $M$ without boundary.

In local coordinates $x=(x^1,\ldots,x^n)$ our operator reads
\begin{equation}
\label{operator in local coordinates}
L=P^\alpha(x)\frac\partial{\partial x^\alpha}+Q(x),
\end{equation}
where $P^\alpha(x)$ and $Q(x)$ are some $m\times m$ matrix-functions
and summation is carried out over $\alpha=1,\ldots,n$.
The full symbol of the operator $L$ is the matrix-function
\begin{equation}
\label{definition of the full symbol}
L(x,p):=iP^\alpha(x)\,p_\alpha+Q(x),
\end{equation}
where $p=(p_1,\ldots,p_n)$ is the dual variable (momentum).

The problem with the full symbol
(\ref{definition of the full symbol})
is that it is not invariant under changes
of local coordinates. The standard analytic way of overcoming this problem is
by introducing the principal and subprincipal symbols in accordance with formulae
\begin{equation}
\label{definition of the principal symbol}
L_\mathrm{prin}(x,p):=iP^\alpha(x)\,p_\alpha\,,
\end{equation}
\begin{equation}
\label{definition of the subprincipal symbol}
L_\mathrm{sub}(x):=Q(x)
+\frac i2(L_\mathrm{prin})_{x^\alpha p_\alpha}(x),
\end{equation}
where the subscripts indicate partial derivatives.
It is known that $L_\mathrm{prin}$ and $L_\mathrm{sub}$
are invariantly defined matrix-functions on $T^*M$ and $M$ respectively,
see subsection 2.1.3
in \cite{mybook} for details.
As we assumed our operator $L$ to be formally self-adjoint,
the matrix-functions $L_\mathrm{prin}$ and $L_\mathrm{sub}$
are Hermitian.

The definition of the subprincipal symbol (\ref{definition of the subprincipal symbol})
originates from the classical paper \cite{DuiHor} of
J.J.~Duistermaat and L.~H\"ormander: see formula (5.2.8) in that paper.
Unlike \cite{DuiHor}, we work with matrix-valued symbols, but this
does not affect the formal definition of the subprincipal symbol.

A peculiar feature of first order differential operators,
as opposed to pseudo\-differential operators and higher order differential operators,
is that the principal and subprincipal symbols uniquely determine the operator.
Namely, examination of formulae
(\ref{operator in local coordinates}),
(\ref{definition of the principal symbol})
and
(\ref{definition of the subprincipal symbol})
gives, in local coordinates, the following expression for the operator in terms of its
principal and subprincipal symbols:
\begin{equation}
\label{operator in terms of its principal and subprincipal symbols}
L=-i[(L_\mathrm{prin})_{p_\alpha}(x)]\frac\partial{\partial x^\alpha}
-\frac i2(L_\mathrm{prin})_{x^\alpha p_\alpha}(x)
+L_\mathrm{sub}(x)\,.
\end{equation}

\section{Massless Dirac action}
\label{Massless Dirac action}

In this appendix we define, in a concise manner, the massless Dirac action.
For more details see Appendix A in \cite{jst_part_b}.

In order to write down the massless Dirac action we need Pauli matrices,
i.e.~a triple of trace-free Hermitian $2\times2$ matrix-functions $\sigma^\alpha(x)$,
$\alpha=1,2,3$, satisfying the condition (\ref{defining relation for Pauli matrices}).
In our case we have Pauli matrices $\sigma^\alpha(x)$ readily available: these are defined in
accordance with formula (\ref{definition of Pauli matrices sigma}),
or, equivalently,
in accordance with formulae
(\ref{definition of Pauli matrices sigma alternative})
and
(\ref{principal symbol of B is linear in xi}).
Covariant Pauli matrices are defined as
$\sigma_\alpha:=g_{\alpha\beta}\sigma^\beta$.

The massless Dirac operator is the matrix operator
\begin{equation}
\label{definition of Weyl operator}
W:=-i\sigma^\alpha
\left(
\frac\partial{\partial x^\alpha}
+\frac14\sigma_\beta
\left(
\frac{\partial\sigma^\beta}{\partial x^\alpha}
+\left\{{{\beta}\atop{\alpha\gamma}}\right\}\sigma^\gamma
\right)
\right),
\end{equation}
where summation is carried out over
$\alpha,\beta,\gamma=1,2,3$, and
\begin{equation}
\label{definition of Christoffel symbols}
\left\{{{\beta}\atop{\alpha\gamma}}\right\}:=
\frac12g^{\beta\delta}
\left(
\frac{\partial g_{\gamma\delta}}{\partial x^\alpha}
+
\frac{\partial g_{\alpha\delta}}{\partial x^\gamma}
-
\frac{\partial g_{\alpha\gamma}}{\partial x^\delta}
\right)
\end{equation}
are the Christoffel symbols.
The operator (\ref{definition of Weyl operator})
acts on a 2-component complex-valued spinor field
$\xi$ which we write as a 2-column,
$\xi=\begin{pmatrix}\xi^1&\xi^2\end{pmatrix}^T$.

We chose the letter ``$W$'' for denoting the massless Dirac operator
because in theoretical physics literature it is often referred to as the \emph{Weyl}
operator. Note that one should really be referring here to the
\emph{static} Weyl operator because we have excluded time,
which is natural in the setting of spectral theory.

We define the massless Dirac action as
\begin{equation}
\label{definition of massless Dirac action}
S(\xi):=\int_M\operatorname{Re}(\xi^*W\xi)\,\sqrt{\det g_{\alpha\beta}}\ dx\,,
\end{equation}
where the star indicates Hermitian conjugation.
This is the variational functional corresponding to the operator
(\ref{definition of Weyl operator}).
Here, of course, we use the fact that in view of the self-adjointness of the operator $W$ we have
\[
\int_M\xi^*(W\xi)\,\sqrt{\det g_{\alpha\beta}}\ dx
=
\int_M(W\xi)^*\xi\,\sqrt{\det g_{\alpha\beta}}\ dx
=
\int_M\operatorname{Re}(\xi^*W\xi)\,\sqrt{\det g_{\alpha\beta}}\ dx\,.
\]

\section{Example}
\label{Example}

In this appendix we consider an explicit example illustrating
the use of  Theorem~\ref{main theorem}.

Consider the unit torus $\mathbb{T}^3$ parameterized by cyclic coordinates $x^\alpha$,
$\alpha=1,2,3$, of period $2\pi$. Let $L$ be the differential operator
with principal symbol
\begin{equation}
\label{Example equation 1}
L_\mathrm{prin}(x,p)=
\begin{pmatrix}
p_3&e^{2ix^3}(p_1-ip_2)\\
e^{-2ix^3}(p_1+ip_2)&-p_3
\end{pmatrix}
\end{equation}
and zero subprincipal symbol.
We examine below the eigenvalue problem
(\ref{eigenvalue problem})
for this particular operator $L$ and trivial weight function
(\ref{case of trivial weight function}).

Substituting
(\ref{Example equation 1})
into
(\ref{definition of metric})
we see that the above principal symbol generates the Euclidean metric
\begin{equation}
\label{Example equation 2}
g^{\alpha\beta}(x)=\delta^{\alpha\beta}.
\end{equation}
Hence, as the reference principal symbol it is natural to take
\begin{equation}
\label{Example equation 3}
\mathring L_\mathrm{prin}(x,p)=
\begin{pmatrix}
p_3&p_1-ip_2\\
p_1+ip_2&-p_3
\end{pmatrix}.
\end{equation}
Substituting
(\ref{Example equation 3})
into
(\ref{definition of Pauli matrices sigma})
we get standard Pauli matrices
\begin{equation}
\label{Example equation 4}
\sigma^1=
\begin{pmatrix}
0&1\\
1&0
\end{pmatrix}
,
\qquad
\sigma^2=
\begin{pmatrix}
0&-i\\
i&0
\end{pmatrix}
,
\qquad
\sigma^3=
\begin{pmatrix}
1&0\\
0&-1
\end{pmatrix}
.
\end{equation}

It is not \emph{a priori} obvious that
the principal symbols
$L_\mathrm{prin}(x,p)$ and $\mathring L_\mathrm{prin}(x,p)$
have the same spin structure. The only way to establish
that they do indeed have the same spin structure
is to resolve equation (\ref{principal symbol via reference principal symbol})
with respect to the special unitary matrix-function $R(x)$.
Straightforward calculations give
\begin{equation}
\label{Example equation 5}
R(x)=\pm
\begin{pmatrix}
e^{ix^3}&0\\
0&e^{-ix^3}
\end{pmatrix}.
\end{equation}
Of course, the underlying reasons why in this particular
case we do not encounter topological
obstructions are that both principal symbols have the same (positive)
topological charge and that the frame encoded in
(\ref{Example equation 1}) makes an even number of turns (two turns)
as $x^3$ runs from 0 to $2\pi$. See Appendix A in \cite{jst_part_b}
for more details.

Formulae
(\ref{SU(2) matrix expressed via spinor}),
(\ref{normalisation of spinor}),
(\ref{case of trivial weight function})
and
(\ref{Example equation 5})
give us the following expression for the spinor field:
\begin{equation}
\label{Example equation 6}
\xi(x)=\pm\begin{pmatrix}e^{-ix^3}\\0\end{pmatrix}.
\end{equation}
Substituting formulae
(\ref{Example equation 2}),
(\ref{Example equation 4})
and
(\ref{Example equation 6})
into
(\ref{definition of Weyl operator})--(\ref{definition of massless Dirac action})
we conclude that $S(\xi)=-(2\pi)^3$.
Hence, Theorem~\ref{main theorem} tells us that in our example
the two-term asymptotics
(\ref{two-term asymptotic formula for counting function})
takes the form
\begin{equation}
\label{Example equation 7}
N(\lambda)=\frac43\pi\lambda^3-4\pi\lambda^2+o(\lambda^2)
\end{equation}
as $\lambda\to+\infty$.
Note that the nonperiodicity condition
(see Definitions 8.3 and 8.4 in \cite{jst_part_a})
is fulfilled in our example, so,
according to Theorem 8.4 from \cite{jst_part_a},
the asymptotic formula
(\ref{Example equation 7})
holds as it is, without mollification.

Observe now that in our example the spectrum of the operator $L$
can be evaluated explicitly. Indeed, let $\mathring L$ be the differential operator
with principal symbol (\ref{Example equation 3})
and zero subprincipal symbol. In other words,
let $\mathring L=\mathring L_\mathrm{prin}(x,-i\partial/\partial x)$.
Consider now the operator
$R\mathring LR^*$, where $R$ is the matrix-function
(\ref{Example equation 5})
It is easy to check that the subprincipal symbol
of the operator $R\mathring LR^*$ is $-I$,
where $I$ is the $2\times2$ identity matrix. Hence,
\begin{equation}
\label{Example equation 8}
L=R\mathring LR^*+I,
\end{equation}
compare with formula (\ref{inequality}).
But the operator $R\mathring LR^*$ is unitarily equivalent to the operator $\mathring L$
and the spectrum of $\mathring L$ is known, see Appendix B in \cite{jst_part_b}.
Using (\ref{Example equation 8}), we conclude that the
eigenvalues of our operator $L$ are as follows.
\begin{itemize}
\item
The number $1$ is an eigenvalue of multiplicity two.
\item
For each $m\in\mathbb{Z}^3\setminus\{0\}$ we have
the eigenvalue $1+\|m\|$ and unique (up to rescaling) eigenfunction,
with eigenfunctions corresponding to different $m$ being linearly
independent.
\item
For each $m\in\mathbb{Z}^3\setminus\{0\}$ we have
the eigenvalue $1-\|m\|$ and unique (up to rescaling) eigenfunction,
with eigenfunctions corresponding to different $m$ being linearly
independent.
\end{itemize}
Thus, $N(\lambda)-1$ is the number of integer lattice points
inside a 2-sphere of radius $\lambda-1$ in $\mathbb{R}^3$ centred at the origin.
Applying the result from \cite{heath-brown_1999} we get
\begin{equation}
\label{Example equation 9}
N(\lambda)=\frac43\pi\lambda^3-4\pi\lambda^2+O_\varepsilon(\lambda^{21/16+\varepsilon})
\end{equation}
as $\lambda\to+\infty$, with $\varepsilon$ being an arbitrary positive number.
The more advanced number theoretic result
(\ref{Example equation 9}) agrees with our
asymptotic formula (\ref{Example equation 7}).

Note that the calculations presented in this section remain
unchanged if we replace everywhere $p_1\mp ip_2$ by
$p_1\pm ip_2$. This is in agreement with the fact that the
topological charge $\mathbf{c}$ does not appear  in our formula
(\ref{formula for b}).

\end{document}